\documentclass[12pt]{report}

\reversemarginpar

\newcommand{\V}{\vskip0.25in}

\begin{document}
\newtheorem{theorem}{Theorem}
\newtheorem{lemma}{Lemma}
\newtheorem{prop}{Proposition}
\newtheorem{cor}{Corollary}
\newtheorem{fact}{Fact}

\begin{center}
{\bf A New Approach to Rational Values of Trigonometric Functions}
\V
{by Greg Dresden, {\tt dresdeng@wlu.edu}, Washington \& Lee University }
\end{center}
\thispagestyle{empty}

\vskip0.5in
\renewcommand{\thefootnote}{\fnsymbol{footnote}}

{\bf NOTE:} a much shorter proof of my Fact \ref{f1}, also using the cyclotomic
polynomials, appears in ``Rational 
Values of Trigonometric Functions\rq\rq\ by Kaoru Motose, MAA
Monthly (114) November 2007, page 818. 

\vskip0.5in

For $a,b$ both integers, when is
$\sin(a\pi/b)$
a rational number?
For that matter, what about $\tan$ and $\cos$? 
We all know about the ``obvious values\rq\rq\ of 
$a$ and $b$ that will give rational answers:
\begin{eqnarray*}
\sin(0) = \tan(0) &=& 0 \qquad \qquad \cos(0)=1 \\
\sin(\pi/6) = \cos(\pi/3) &=& 1/2  \\
\tan(\pi/4) &=& 1\\
\cos(\pi/2) &=& 0 \qquad\qquad \sin(\pi/2) = 1
\end{eqnarray*}

\noindent ...and so on. (For ease of discussion, let's keep $a/b$
in the interval $[0, 1/2]$.)

Are there any \underbar{other} values for $a/b$ such that
$\sin(a\pi/b)$
(or $\cos$, or $\tan$) is rational? 
The answer, of course, is no
(as one colleague quipped, if it was rational anywhere else, we
surely would have heard about it!).
Let us express this fact in the following informal manner:
\begin{fact}\label{f1}
For $a,b$ relatively prime integers (with $b > 0$), 
then 
$\sin(a\pi/b)$, 
$\cos(a\pi/b)$, and 
$\tan(a\pi/b)$ are rational only at the
obvious values of $a/b$ (in particular, $b$ can not be other than $1,2,3,4$, or $6$).
\end{fact}
The classical proofs of this fact involve
the Chebyshev polynomials and various 
trig 
identities (see [1], [3, section 6.3], and
[5], as well as 
the commentary after [4]). Chebyshev polynomials rarely appear in 
the traditional
undergraduate curriculum, and thus the proof of Fact \ref{f1} is
not usually seen by students. 
In this paper, we utilize a different procedure, and show that Fact \ref{f1}
is in fact equivalent to the following well-known statement, 
familiar to most algebra students:
\begin{fact}\label{f2}
For $c,d$ relatively prime integers (with $d > 0$), the primitive
$d$th root of unity $e^{2\pi ic/d}$  has degree $\leq 2$ over
$\mathbf{Q}$
iff
$d=1,2,3,4$, or $6$.
\end{fact}

We point out that this topic 
is well suited for an abstract algebra class, and
provides a delightful application of the theory of field
extensions. 
The method outlined here 
is relatively straightforward and 
would
involve only a few minutes of classroom time (alternatively, it 
would make an excellent homework assignment).
Indeed, proving Fact \ref{f2} independently takes very little
work; 
one might first show that
the degree of the field 
$\mathbf{Q}(e^{2\pi ic/d})$ is 
$\phi(d)$ (perhaps by showing that the cyclotomic polynomial
$\Phi_d(x)$ is irreducible) and one could then
show that $\phi(d) \leq 2$ only for the values
of $d$ given above. We leave the details as an
exercise for the reader (see [2, chapter 33]).

\V

Let us now show the equivalence of our two facts.

First, suppose Fact \ref{f1} is true. Let $c,d$ 
be relatively prime integers (with $d>0$), and suppose
$K=\mathbf{Q}(e^{2\pi ic/d})$ is  of degree $2$ over $\mathbf{Q}$.
By Euler\rq s formula, we can write this primitive $d$th
root of unity as  
$e^{2\pi ic/d} = \cos(2\pi c/d) + i\sin(2\pi c/d)$.
With this in mind, we note that
the field $K$ contains the real number
$(1/2)\left(e^{2\pi ic/d} + 1/e^{2\pi ic/d}\right) = \cos(2\pi c/d)$
and thus also the complex number $i\sin(2\pi c/d)$. These can\rq t
both be degree $2$ over $\mathbf{Q}$, as 
the field $K$, being only of degree $2$,
 can\rq t contain 
both a
real degree-$2$ subfield
and a 
complex degree-$2$ subfield.
Thus, either 
$\sin(2\pi c/d)=0$ or
$\cos(2\pi c/d) \in \mathbf{Q}$. By 
Fact \ref{f1}, the first case gives
$d=1$ or $2$, and the second gives 
$d=1,2,3,4$, or $6$.

\V
Second, suppose Fact \ref{f2} is true.  Choose a rational number $a/b$ in
reduced form such that $\tan(a\pi/b)$ equals some
rational number $r$, and
let $v=1+ri$ (see Figure 1, below). 
Now, $v$ is in $\mathbf{Q}(i)$, but since it\rq s not of 
length $1$, it clearly is not a root of unity and so we can\rq t use
Fact \ref{f2}.
So, it would be reasonable to consider
\[
\frac{v}{|v|} = \frac{1}{\sqrt{1+r^2}} + 
\frac{r}{\sqrt{1+r^2}}i = e^{\pi ia/b},
\]
which clearly has length $1$ and argument 
$a\pi/b$, and thus is a root of unity.
Unfortunately, this complex number is in the possibly
degree-$4$ field $\mathbf{Q}(\sqrt{1+r^2},i)$ so 
we still can\rq t apply Fact \ref{f2}! Instead, we look at
\[
\left(\frac{v}{|v|}\right)^2 = 
\frac{1-r^2}{{1+r^2}} + \frac{2r}{{1+r^2}}i = e^{2\pi ia/b},
\]
which is clearly in the quadratic number field $\mathbf{Q}(i)$. Thus, by
Fact \ref{f2} (and since $a,b$ are relatively prime)
we have that $b = 1,2,3,4$, or $6$; a simple
calculation shows that tangent is rational only at the
obvious values.

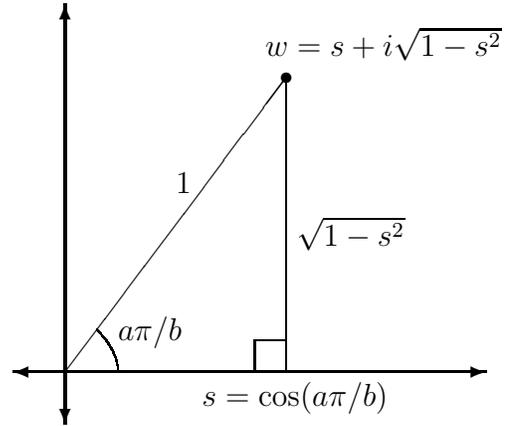
\begin{minipage}{90mm}
\setlength{\unitlength}{1.4mm}
\begin{picture}(50,50)
\thicklines
\put(5,5){\vector(1,0){40}}
\put(5,5){\vector(0,1){35}}
\put(5,5){\vector(-1,0){5}}
\put(5,5){\vector(0,-1){5}}
\thinlines
\put(26,33){\circle*{1}}
\put(24,35){$v = 1+ri$}
\put(27,17){$r = \tan(a\pi/b)$}
\put(17,22){\makebox(0,0)[br]{$\sqrt{1+r^2}$}}
\put(18,4){\makebox(0,0)[tl]{$1$}}
\put(25,43){\makebox(0,0)[bc]{\bf Figure 1}}
\put(10,8){$a\pi/b$}
\put(5,5){\line(3,4){21}}
\put(26,5){\line(0,1){28}}
\put(23,5){\line(0,1){3}}
\put(26,8){\line(-1,0){3}}
\qbezier(10,5)(10,7)(8,9)
\end{picture}
\end{minipage}
\hfill 
\begin{minipage}{90mm}
\setlength{\unitlength}{1.4mm}
\begin{picture}(50,50)
\thicklines
\put(5,5){\vector(1,0){40}}
\put(5,5){\vector(0,1){35}}
\put(5,5){\vector(-1,0){5}}
\put(5,5){\vector(0,-1){5}}
\thinlines
\put(26,33){\circle*{1}}
\put(25,43){\makebox(0,0)[bc]{\bf Figure 2}}
\put(24,35){$w = s+i\sqrt{1-s^2}$}
\put(27,17){$\sqrt{1-s^2}$}
\put(18,4){\makebox(0,0)[tl]{$s = \cos(a\pi/b)$}}
\put(10,8){$a\pi/b$}
\put(17,22){\makebox(0,0)[br]{$1$}}
\put(5,5){\line(3,4){21}}
\put(26,5){\line(0,1){28}}
\put(23,5){\line(0,1){3}}
\put(26,8){\line(-1,0){3}}
\qbezier(10,5)(10,7)(8,9)
\end{picture}
\end{minipage}

\vskip0.3in

We now proceed to show the same holds for 
cosine (once we have this, the rationality of
sine follows from the identity
$\sin(\theta) = \cos(\pi/2 - \theta)$.
In a similar manner to our work earlier, we 
choose a rational number $a/b$ in
reduced form such that $\cos(a\pi/b) = s$ 
(for $s$ some rational number), and let
$w=s+i\sqrt{1-s^2}$ (see Figure 2, above).
 Now, 
$|w| = 1$ and 
$\arg(w) = a\pi/b$, so 
$w = e^{i a\pi/b}$ and is in $\mathbf{Q}(i\sqrt{1-s^2})$,
a (complex) quadratic number field. 
Thus, we can apply Fact \ref{f2} to note that $b$ must be
$1,2,3,4$, or $6$, and again, calculations give us the 
desired obvious values.

\V
This completes our proof
of the equivalence of the two facts, but it does not mark the
end of this intriguing area of study. For example, we note that the
roots of unity of degree $4$ are those numbers
$e^{2\pi ic/d}$ with 
$d \in \phi^{-1}(4) = \{5,8,10,12\}$. 
Likewise, we note that:
\begin{eqnarray*}
\cos(\pi/5) = \frac{\sqrt{5}+1}{4} &\qquad \qquad &  
\sin(\pi/10) = \frac{\sqrt{5}-1}{4}\\
\tan(\pi/8) = \sqrt{2}-1 
&\qquad \qquad &  
\tan(\pi/12) = 2-\sqrt{3},
\end{eqnarray*}
all simple radicals of degree $2$ over
$\mathbf{Q}$. 
The interested reader might want to generalize Facts 1 and 2
to include this correspondence (as well as others of
arbitrary degree).
Indeed, this might well lead to an alternate proof
of the well-known statement that
the trig functions are algebraic at all rational multiples 
of $\pi$.

\vskip0.5in

\vskip1.0in

\noindent {\sc references}

\begin{enumerate}

\item L. Carlitz and J. M. Thomas, 
\emph{Rational Tabulated Values 
of Trigonometric Functions}, 
{Amer. Math. Monthly} {\bf 69} (1962) pp. 789--793.  

\item J. Gallian, \emph{Contemporary Abstract Algebra} (5th ed.), 
Houghton-Mifflin, Boston, 2002.

\item I. Niven, H. Zuckerman, and H. Montgomery,
\emph{An Introduction to the Theory of Numbers} (5th ed.),
Wiley \& Sons, New York, 1991.

\item  Andrzej Makowski, D. E. Penney, 
\emph{Problem E2333},
Amer. Math. Monthly {\bf 80} (1973) pp. 77-78.

\item J. M. H. Olmsted, 
\emph{Rational Values 
of Trigonometric Functions}, 
 Amer. Math. Monthly {\bf 52} (1945) pp. 507--508.

\end{enumerate}
\end{document}